\documentclass[letterpaper, paper,11pt]{AAS}		
\usepackage{tcolorbox}
\usepackage{lipsum} 
\usepackage{bm}
\usepackage{amsmath,amssymb}
\usepackage{float}
\usepackage[colorlinks=true, pdfstartview=FitV, linkcolor=black, citecolor= black, urlcolor= black]{hyperref}
\usepackage{subcaption}
\usepackage{overcite}
\usepackage{footnpag}			      	

\PaperNumber{23-292}

\begin{document}

\author{Simone Semeraro\thanks{Department of Aerospace Engineering, The University of Michigan, {Ann Arbor}, {Michigan}, {USA}},  
Ilya Kolmanovsky\thanks{Corresponding Author. Department of Aerospace Engineering, The University of Michigan, {Ann Arbor}, {Michigan}, {USA}. {This research was supported by National Science Foundation Grant Numbers CMMI-1904394 and ECCS-1931738.}}, Emanuele Garone\thanks{École Polytechnique de Bruxelle, {Brussels}, Belgium}}

\title{ On Constrained Feedback Control of Spacecraft Orbital Transfer Maneuvers}


\maketitle{}

\begin{abstract}

The paper revisits a Lyapunov-based feedback control to implement spacecraft orbital transfer maneuvers.  The spacecraft equations of motion in the form of Gauss Variational Equations (GVEs) are used.  By shaping the Lyapunov function using barrier functions, we demonstrate that state and control constraints during orbital maneuvers can be enforced.  Simulation results from orbital maneuvering scenarios are reported. The synergistic use of the reference governor in conjunction with the barrier functions is proposed to ensure convergence to the target orbit (liveness) while satisfying the imposed constraints. 

\end{abstract}

\begin{tcolorbox}[colback=gray!10!white, colframe=gray!50!black, width=\textwidth, boxrule=0.5mm, sharp corners]
	\begin{center}
		\small
		{This document is an author's original manuscript (pre-print) of the paper AAS 23-292 presented at the 33rd AAS/AIAA Space Flight Mechanics Meeting, Austin, Texas, January 15-19, 2023.} \\
		\scriptsize
		\url{https://www.space-flight.org/docs/2023_winter/33rdSFFM-Full_Program_Abstracts.pdf}
		\normalsize
	\end{center}
\end{tcolorbox}

\section{Introduction}
This article considers the development of feedback laws based on Gauss Variational Equations (GVEs) for stabilization of the spacecraft to the target orbit subject to constraints on minimum radius of periapsis, thrust magnitude and orbit eccentricity.
The spacecraft is assumed to have continuous thrust capability. In the case of on-off thrusters, continuous thrust values are assumed to be realized using pulse width modulation\cite{wie1998space}.

The feedback control for orbital transfers has  been previously investigated \cite{gurfil2016celestial,petropoulos2005refinements,holt2020low,chang2002lyapunov} for missions that exploit low thrust propulsion. In particular, the Q-law \cite{petropoulos2004low,petropoulos2005refinements} has been  extensively studied.
The feedback control laws provide corrections all the time as opposed to only at pre-defined Trajectory Correction Maneuver (TCM) points. The use of feedback for orbital transfers is appealing due to its potential to improve robustness to unmodeled perturbations and thrust errors.  At the same time, nonlinear dynamics of spacecraft motion and constraints, such as thrust limits that are much smaller as compared to the dominant gravitational forces, complicate the application of nonlinear control methods.

Our approach to enforce the constraints involves modifying the nominal Lyapunov feedback law in
the unconstrained case\cite{gurfil2016celestial} with the barrier functions.  We also propose the use of the reference governor\cite{garone2017reference} to ensure the convergence of the trajectory to the target orbit.

\section{Modeling}

The Gauss-Euler Variational Equations (GVEs) \cite{gurfil2016celestial} are exploited as equations of motions that apply to the six classical orbital elements of Keplerian two body problem in presence of perturbation forces such as spacecraft thrust.
These orbital elements are  the semi-major axis $a$ [km], the eccentricity $e$,  the inclination $i$ [rad],  the Right Ascension of the Ascending Node (RAAN) $\Omega$ [rad], the argument of periapsis $\omega$ [rad] and the spacecraft true anomaly $\theta$ [rad].  

To model the evolution of the orbital elements, a moving STW frame is introduced with the origin at the spacecraft's center of mass and with the unit vectors $\hat{e}_r$, $\hat{e}_\theta$ and $\hat{e}_h$
defined according to
$$\hat{e}_r=\frac{\vec{r}}{r},~\hat{e}_h=\frac{\vec{h}}{h},~ \vec{h} =\vec{r} \times \vec{v},~
\mbox{and}~ \hat{e}_\theta=\hat{e}_h \times \hat{e}_r,$$
where $\vec{r}$ is the spacecraft position vector in an inertial frame centered at the primary body,  $\vec{v}$ is the spacecraft velocity, and  $\times$ denotes the vector product.
Decomposing the thrust force per unit mass 
[km/s$^2$]  applied to the spacecraft as
\begin{align*}
    \frac{\vec{F}}{m} = S\hat{e}_r + T\hat{e}_\theta + W\hat{e}_h,
\end{align*}
where $m$ denotes the mass of the spacecraft,
the GVEs take the following form:
\begin{align}\label{gve}
    \frac{da}{dt} &= \frac{2a^2}{\sqrt{\mu p}}e\sin{\theta}S+\frac{2a^2}{\sqrt{\mu p}}\frac{p}{r}T, \nonumber \\
    \frac{de}{dt} &= \frac{p\sin\theta}{\sqrt{\mu p}}S + \frac{p(\cos\psi+\cos\theta)}{\sqrt{\mu p}}T, \nonumber\\
    \frac{di}{dt} &= \frac{r}{\sqrt{\mu p}}\cos{(\theta+\omega)}W, \\
    \frac{d\Omega}{dt} &= \frac{r}{\sqrt{\mu p}}\frac{\sin{(\theta+\omega)}}{\sin i}W, \nonumber\\
    \frac{d\omega}{dt} &= -\frac{p\cos\theta}{e\sqrt{\mu p}}S+\frac{(r+p)\sin\theta}{e\sqrt{\mu p}}T-\frac{r\sin(\theta+\omega)\cot i}{\sqrt{\mu p}}W,\nonumber \\
    \frac{d\theta}{dt} &= \frac{\sqrt{\mu p}}{r^2}+\frac{p\cos\theta}{e}\frac{S}{\sqrt{\mu p}}-\frac{p+r}{e}\cos\theta\frac{T}{\sqrt{\mu p}}, \nonumber
\end{align}
where $\mu$ is the gravitational parameter of the primary, $r=\frac{p}{1+e\cos\theta}$ the distance from the gravity center to the spacecraft center of mass, $p=a(1-e^2)$ the orbit parameter (also called semi-latus rectum), and $\psi = \arccos\left(\frac{1}{e}-\frac{r}{ae}\right)$ is the eccentric anomaly.

Let
\begin{align*}
    X = [a \ e \ i \ \Omega \ \omega]^{\sf T}, \ \ U = [S \ T \ W]^{\sf T}.
\end{align*}
Then the first five of GVEs (\ref{gve}) can be written in the following condensed form:
\begin{equation} \label{equ:main1}
    \dot X(t) = G\big(X(t),\theta(t)\big)U(t),
\end{equation}
where $G(X(t),\theta(t)) \in \mathbb{R}^{5 \times 3}$.
Note that since the objective is to steer the spacecraft into a target orbit but not to a particular location in that orbit (each location will be visited due to periodicity of the orbit) the true anomaly $\theta(t)$ is not included into the state vector $X(t)$ and is treated as a time-varying parameter the evolution of which is determined by the sixth equation
in (\ref{gve}).  Note that the system (\ref{equ:main1}) is an affine nonlinear control system which is drift-free.

\section{Constraints}

Three constraints are considered on the spacecraft trajectory during the orbital transfer.
The first constraint has the form,
\begin{equation}\label{equ:maincon1}
    c_1(X) = r_{\tt p} - r_{\tt min} \geq 0,\quad r_{\tt p}=a(1-e). 
\end{equation}
It ensures that the radius of the periapsis of the spacecraft orbit, $r_{\tt p}$,
is larger than $r_{\tt min}$.  The constraint (\ref{equ:maincon1})
protects not only against the distance to the primary $r$ falling below $r_{\tt min}$ at any given time instant (and thus avoiding being too close to/colliding with the primary) but also that $r$ will stay above $r_{\tt min}$ even if there is a thruster failure and thrust becomes zero.

The second constraint is imposed on the spacecraft thrust magnitude not to exceed the value that the thruster can deliver. Assuming that a feedback law of the form,
$$ U = U(X, X_{\tt des}, \theta),$$
is employed, where $X_{\tt des}$ is the vector of the five orbital elements of the target orbit (target state),
this constraint has the following form,
\begin{equation}\label{equ:maincon2}
    c_2(X, X_{\tt des}, \theta) = U_{\tt max}^2 - \| U \|_2^2 \geq 0,
\end{equation}
ensuring that the spacecraft relative acceleration due to thrust, $U$, remains below a specified value, $U_{\tt max}$, in magnitude.
Here $\|\cdot\|_2$ denotes the standard Euclidean 2-norm; in the sequel $\|\cdot\|=\|\cdot\|_2$ unless specified otherwise.

This constraint assumes that the spacecraft has a single orbital maneuvering thruster and that the attitude of the spacecraft is changed by a faster attitude control loop to accurately realize the commanded thrust direction.  

The design of a Lyapunov-based controller in the case of an $\infty$-norm constraint on $U$ (i.e., if each control channel is constrained individually) and even more general constraints on $U$ is considered in Appendix.  

The third constraint ensures that the eccentricity of the spacecraft orbit is maintained above a specified minimum value, $e \geq e_{\tt min}$, where $e_{\tt min} \geq 0$ :
\begin{equation}\label{equ:maincon3}
    c_3(X) = e-e_{\tt min} \geq  0.
\end{equation}
This constraint ensures that the eccentricity does not approach zero too closely where GVEs have a singularity therefore preserving the validity of the model.  

\section{Control Design}

Given 
the vector of five orbital elements of the target orbit, $X_{\tt des}$, to achieve target tracking while satisfying the constraints (\ref{equ:maincon1}) and (\ref{equ:maincon3}), we
let $P=P^{\sf T} \succ 0$ be a $5 \times 5$ be a positive-definite weight matrix and select a Lyapunov function candidate,
\begin{equation}\label{equ:main2}
    V(X)=\frac{1}{2}\left(X-X_{\tt des}\right)^{\sf T}P\left(X-X_{\tt des}\right)+B_1(X)+B_2(X),
\end{equation}
where
\begin{equation}\label{equ:main2-1}
    B_1(X) = \left\{
    \begin{array}{cc}
         \frac{1}{2} q_1( a(1-e) - r_{\tt min} - \epsilon_1)^2 
         \hfill & 
         \text{if } a(1-e) < r_{\tt min}+\epsilon_1, \\
         0 \hfill & \text{otherwise},
    \end{array} \right. 
\end{equation}

\begin{equation}\label{equ:main2-2}
B_2(X)= \left\{ {\begin{array}{cc}
         \frac{1}{2}q_2(e-e_{\tt min}-\epsilon_2)^2 \hfill & \text{if }e < e_{\tt min}+\epsilon_2, \\
         0 \hfill & \text{otherwise}.
    \end{array}}\right.
\end{equation}
Here $q_1>0$ and $q_2>0$ are weight parameters penalizing the respective constraint violations and $\epsilon_1 \geq 0$, $\epsilon_2 \geq 0$ are safety margins
so that the barrier functions become non-zero prior
to constraints becoming violated.

Computing the time derivative of $V(X)$ along the trajectories of (\ref{equ:main1}) 
we obtain,
\begin{gather*}
    \frac{d}{dt} V(X(t)) = (X(t)-X_{des})^{\sf T} PG(X(t),\theta(t))U(t) + C^{\sf T} (X(t))G(X(t),\theta(t)) U(t),
    \end{gather*}
    where
    \begin{gather*}
    C(X) = \nabla B_1(X) + \nabla B_2(X), \\ \\ \nabla B_1(X)=\\\left\{ {\begin{array}{cc}
         \big[q_1(1-e)(a(1-e)-r_{\tt min}-\epsilon_1), \ -q_1a(a(1-e)-r_{\tt min}-\epsilon_1), \ 0, \ 0, \ 0\big]^{\sf T} \hfill & \text{if }a(1-e) < r_{\tt min}+\epsilon_1, \\
         0_{5\times1}, \hfill & \text{otherwise}, \\
    \end{array}}\right. \\
    \nabla B_2(X)=\left\{ {\begin{array}{cc}
         \big[0, \  q_2(e-e_{\tt min}-\epsilon_2), \ 0, \ 0, \ 0\big]^{\sf T} \hfill & \text{if } e < e_{\tt min}+\epsilon_2, \\
         0_{5\times1}, \hfill & \text{otherwise}. \\
    \end{array}}\right.
\end{gather*}
A Lyapunov control law is now defined as
\begin{equation} \label{equ:controllaw2}
 U = \left\{\begin{array}{ll} U_{\tt nom}, & \mbox{if $\|U_{\tt nom}\| \leq U_{\tt max}$}, \\
 \frac{U_{\tt nom}}{\|U_{\tt nom}\|} U_{\tt max}, & \mbox{otherwise},
 \end{array} \right.
 \end{equation}
 where
\begin{align} \label{equ:controllaw1}
    U_{\tt nom}(X,X_{\tt des},\theta)= -\left[(X-X_{\tt des})^{\sf T} PG(X,\theta) + C^{\sf T}(X)G(X,\theta)\right]^{\sf T}.
\end{align}
With this feedback law,
$$\frac{d}{dt} V(X(t)) \leq 0,$$
along the closed-loop system trajectories implying that the sublevel sets of $V$ are positively invariant.  Furthermore, the constraint (\ref{equ:maincon2}) is enforced due to saturation used in (\ref{equ:controllaw2}).

Consider now a trajectory $X(t)$ evolving from a state, $X(t_0)$, that does not violate the constraints plus safety margins, so that $B_1(X(t_0)) = 0$, $B_2(X(t_0)) = 0$. 
Let 
\begin{equation}\label{equ:V_0}
 V_0(P,X(t_0),X_{\tt des})=\frac{1}{2}\left(X(t_0)-X_{\tt des}\right)^{\sf T}P\left(X(t_0)-X_{\tt des}\right),
 \end{equation}
and note that $V(X(t_0))=V_0(P,X(t_0),X_{\tt des})$.
Since $V(X(t))$ is a non-increasing function and $B_1(X(t))$ and $B_2(X(t))$ are nonnegative,
\begin{align*}
\max &\{B_1(X(t)), B_2(X(t))\} \\ &\leq B_1(X(t))+B_2(X(t)) \leq V(X(t)) \leq V(X(t_0)) \leq  V_0(P,X(t_0),X_{\tt des}).
\end{align*}
If $P$, $q_1$, $q_2$, $\epsilon_1$, $\epsilon_2$ satisfy
\begin{equation}\label{equ:logic}
q_1 \geq \frac{2}{ \epsilon_1^2}V_0(P,X(t_0),X_{\tt des}),~ q_2 \geq \frac{2}{ \epsilon_2^2}V_0(P,X(t_0),X_{\tt des}),
\end{equation}
then $B_1(X(t))\leq \frac{1}{2}q_1 \epsilon_1^2$ and $B_2(X(t)) \leq \frac{1}{2} q_2 \epsilon_2^2$ for all $t \geq t_0$; hence,
the constraints (\ref{equ:maincon1}) and (\ref{equ:maincon3}) are enforced along the closed-loop trajectory.  

The feedback law (\ref{equ:controllaw2})-(\ref{equ:controllaw1}) can be interpreted \cite{kolmanovsky2008speed} as an approximate solution to an optimization problem of minimizing $V(X(t+\Delta t))$, for small $\Delta t >0$, weighted by the control effort, i.e., of
$$J=V(X(t+\Delta t))+\frac{1}{2} \int_t^{t+\Delta t} U^{\sf T}(\tau) U(\tau) d \tau $$
subject to
$$    \|U(\tau)\| \leq U_{\tt max},~t \leq \tau \leq t+\Delta t. $$
This one step ahead predictive control interpretation facilitates tuning the weights $P$, $q_1$ and $q_2$ in (\ref{equ:main2})-(\ref{equ:main2-2}).  For instance, if faster response of one of the orbital elements is desired, the corresponding
diagonal entry of $P$  can be increased.  This process is not dissimilar to the one used for 
Linear Quadratic Regulator (LQR) tuning. A similar procedure of augmenting a 
nominal Lyapunov function with the barrier functions to protect against constraint violations has been previously exploited for other
applications\cite{kolmanovsky2002speed}. The weights $q_1$ and $q_2$ are increased to avoid constraint violations.

We note that conditions (\ref{equ:logic}) can be conservative in that they result in large values of $q_1$ and $q_2$ if set at the beginning of the maneuver. Instead, these conditions  
can be used to reset $q_1$ and $q_2$ periodically
along the trajectory.  Specifically if $B_1(X(t_k))=B_2(X(t_k))=0$ at any time instant $t_k$,
then resetting
\begin{equation}\label{equ:reset}
q_1(t_k) = \frac{2}{ \epsilon_1^2}V_0(P,X(t_k),X_{\tt des}),~ q_2(t_k) =\frac{2}{ \epsilon_2^2}V_0(P,X(t_k),X_{\tt des}),
\end{equation}
does not change the value of the Lyapunov function, $V(X(t_k))$, and still ensures that the constraints are enforced.  As $V_0(P,X(t),X_{\tt des})$ is typically decreasing with time $t$, the periodic reset (\ref{equ:reset})  lowers values of $q_1$ and $q_2$ and facilitates recovering the unconstrained controller performance.  Note also that the values of $q_1$ and $q_2$ will need to be reset whenever $X_{\tt des}$ changes.

\section{Numerical Results}
\begin{figure}[htpb]
\centering
\subfloat[]{
    \includegraphics[width=4cm]{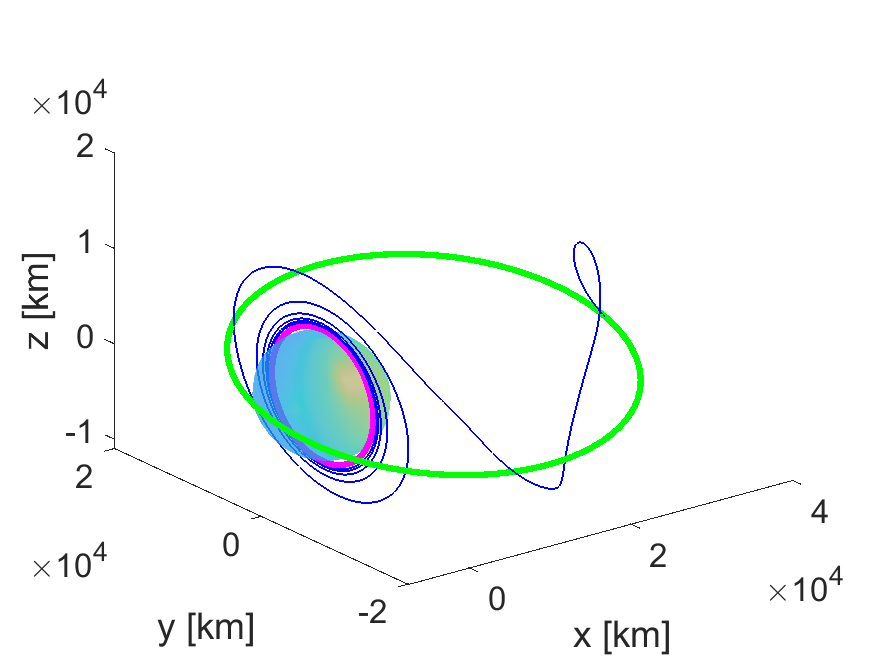}
}\qquad
\subfloat[]{
    \includegraphics[width=4cm]{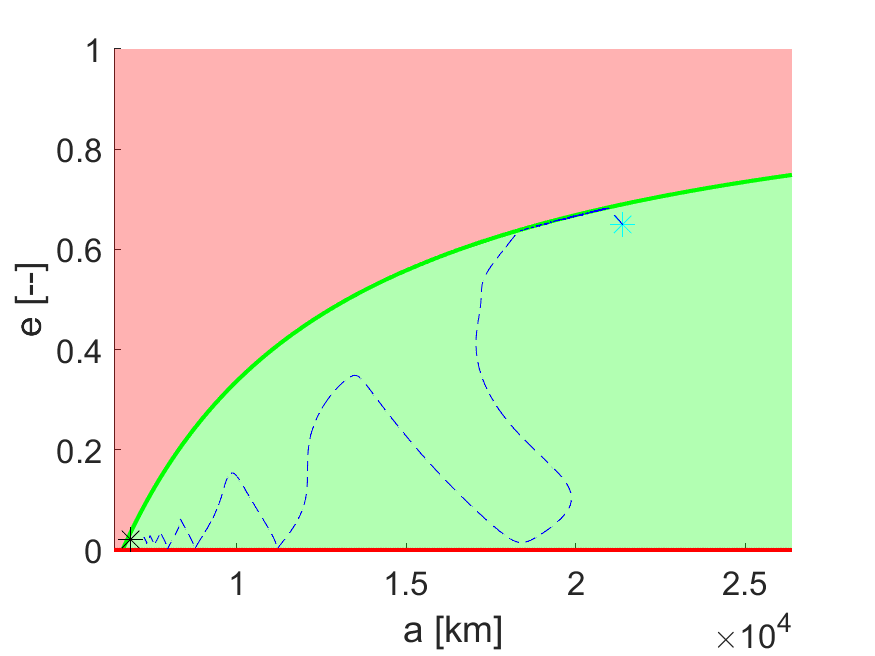}
}\qquad
\subfloat[]{
    \includegraphics[width=4cm]{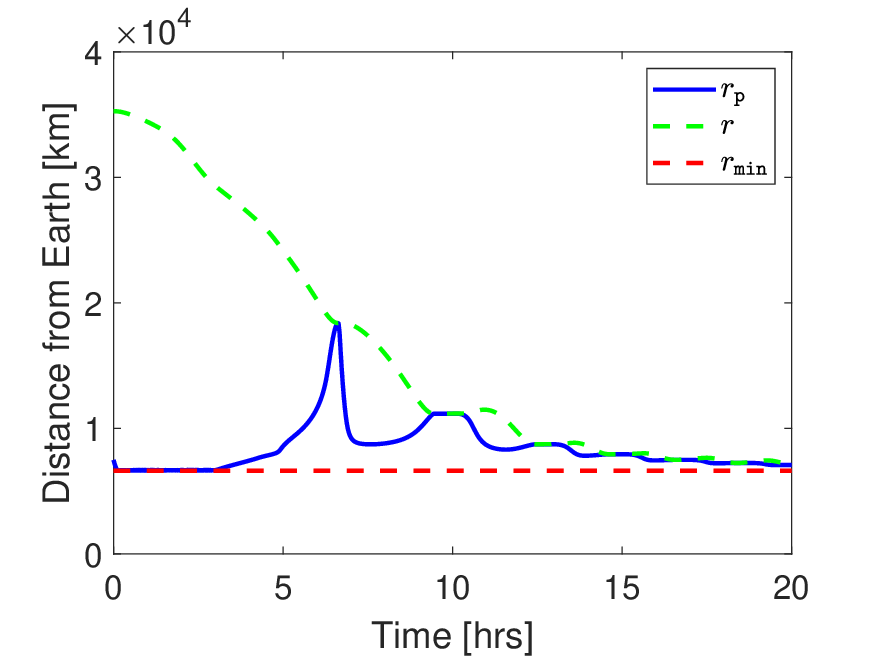}
}\vskip\baselineskip
\subfloat[]{
    \includegraphics[width=4cm]{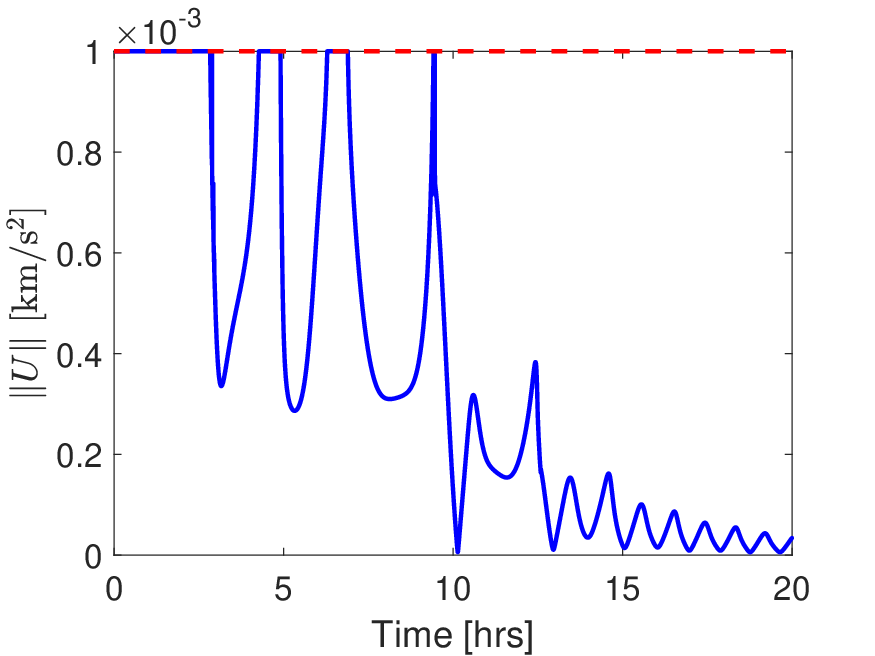}
}\qquad
\subfloat[]{
    \includegraphics[width=4cm]{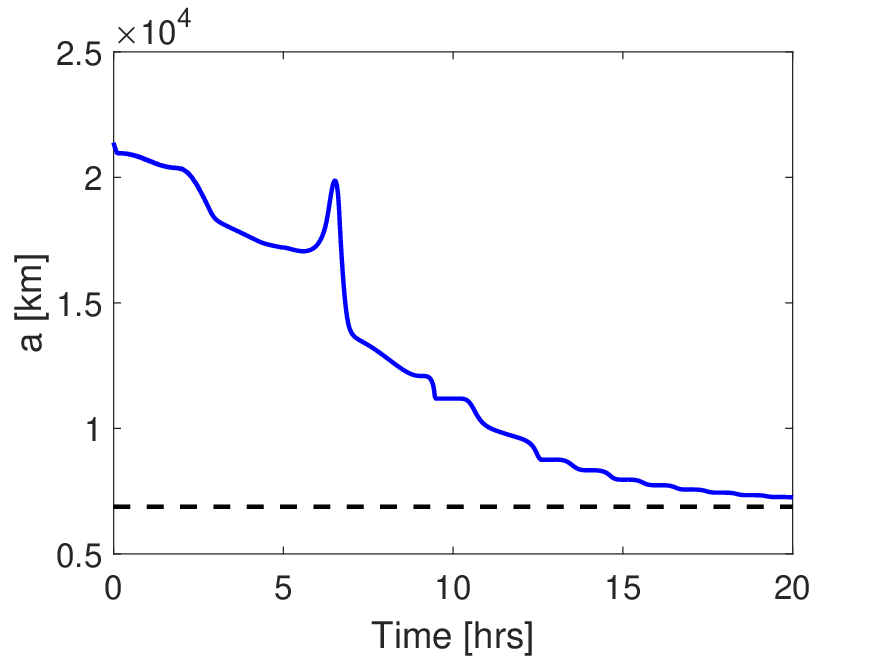}
}\qquad
\subfloat[]{
    \includegraphics[width=4cm]{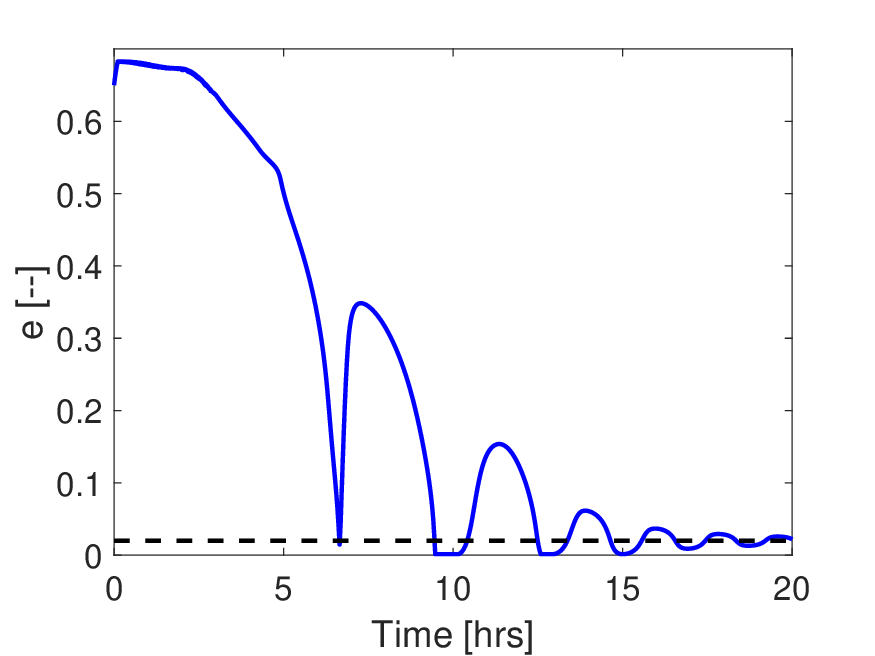}
}\qquad
\subfloat[]{
    \includegraphics[width=4cm]{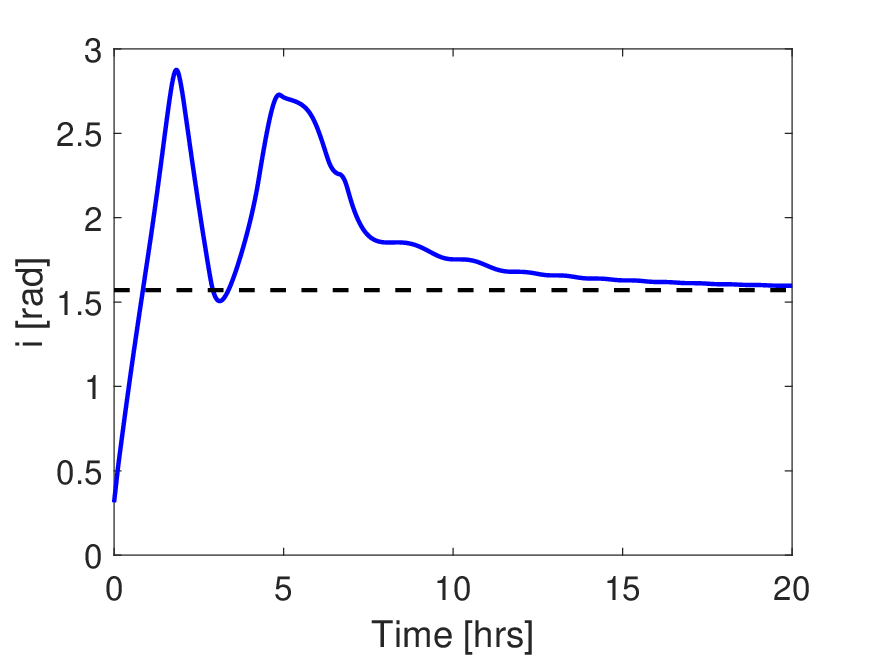}
}\qquad
\subfloat[]{
    \includegraphics[width=4cm]{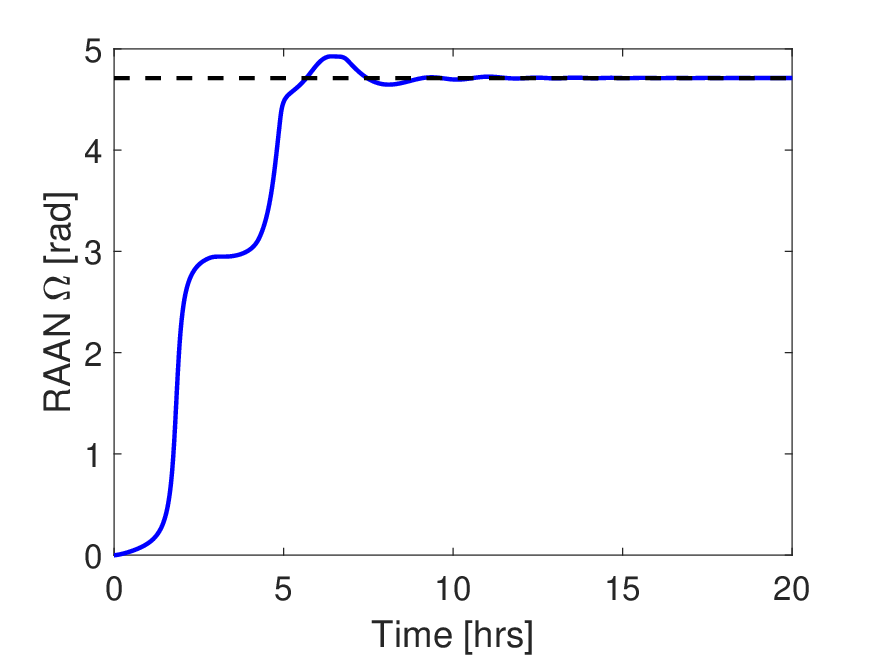}
}\qquad    
\subfloat[]{
    \includegraphics[width=4cm]{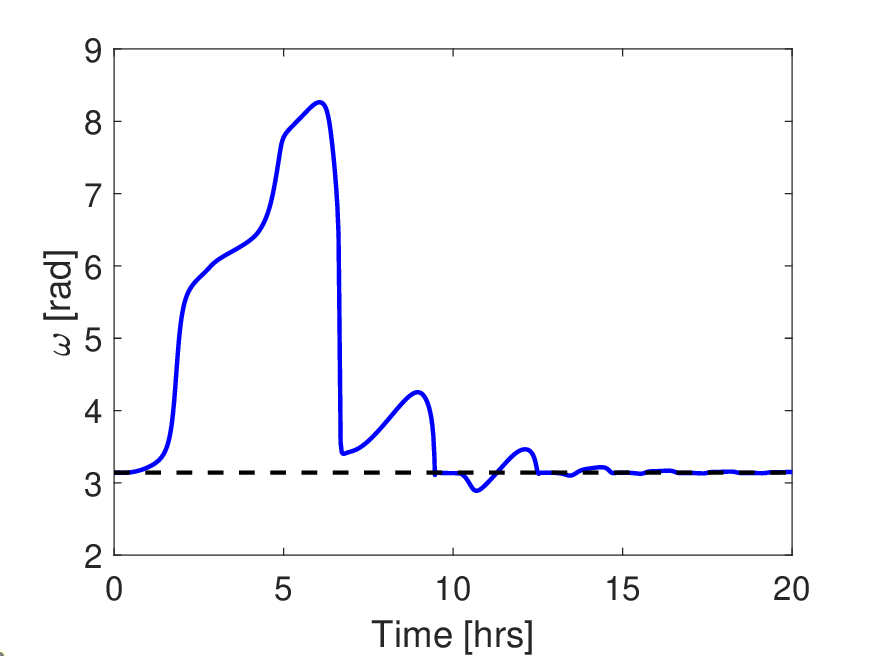}
}\qquad
\subfloat[]{
    \includegraphics[width=4cm]{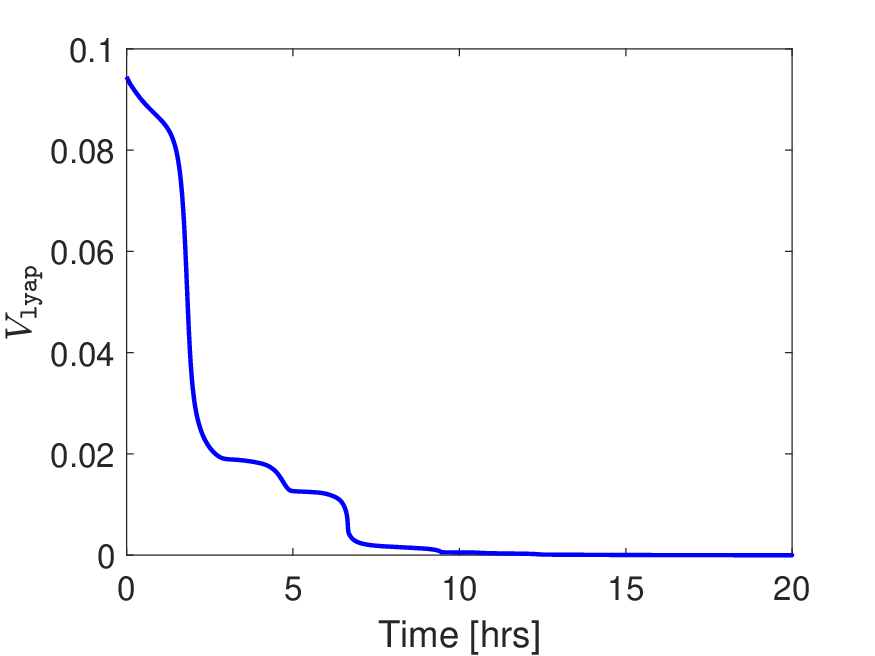}
}\qquad
\subfloat[]{
    \includegraphics[width=4cm]{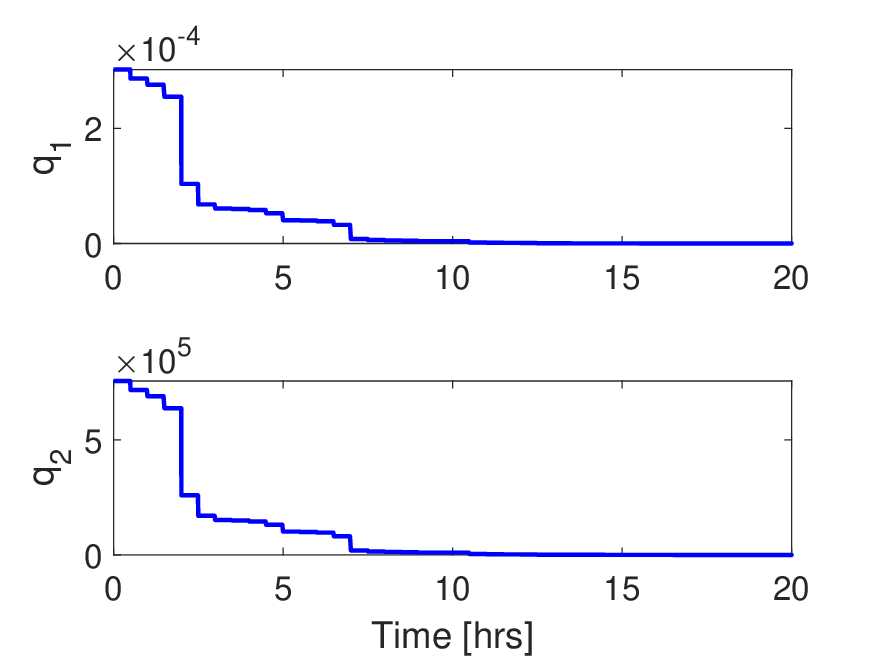}
}\qquad
\subfloat[]{
    \includegraphics[width=4cm]{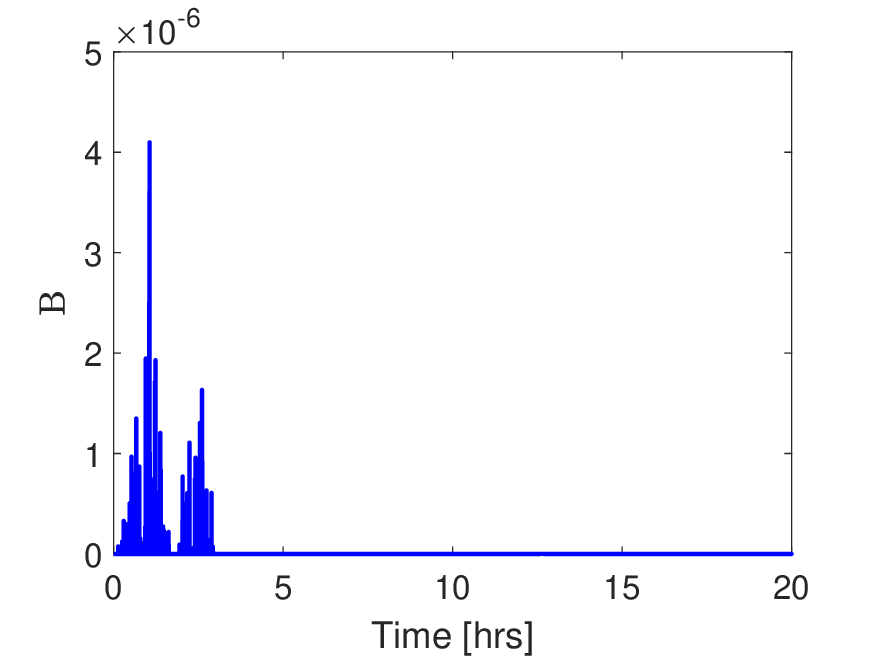}
}
\caption{Orbital transfer from a higher orbit to a lower orbit with an aggressive controller: (a) Three dimensional trajectory; (b) Trajectory on $a$-$e$ plane (dashed) with the region allowed by constraints (\ref{equ:maincon1}) and (\ref{equ:maincon3}) shown in green; (c) The time histories of $r_{\tt p}$, $r$ and $r_{\tt min}$; (d) The time histories of $\|U\|$ (solid)
and $U_{\tt max}$ (dashed); (e) The time histories of $a$ (solid) and $a_{\tt des}$ (dashed); (f) The time histories of $e$ (solid) and $e_{\tt des}$ (dashed); 
(g) The time histories of $i$ (solid) and $i_{\tt des}$ (dashed);
(h) The time histories of $\Omega$ (solid) and $\Omega_{\tt des}$ (dashed);
(i) The time histories of $\omega$ (solid) and $\omega_{\tt des}$ (dashed); (j) The time history of the Lyapunov function, $V(X(t))$;
(k) The time histories of $q_1$ and $q_2$; (l) The time history of the barrier function, $B(X(t))$.}\label{fig:1}
\end{figure}

The maneuver from a higher orbit to a lower  orbit was simulated corresponding to the following initial conditions,
\begin{equation}\label{equ:X0}
    X(0) = \left[21378, \ 0.65, \ \frac{\pi}{10}, \ 0, \ \pi \right]^{\sf T},~~\theta(0)=\pi,
\end{equation}
and desired final conditions corresponding to the target orbit,
\begin{equation}\label{equ:Xdes}
    X_{\tt des} = \left[6878, \ 0.02, \ \frac{\pi}{2}, \ \frac{3\pi}{2}, \ \pi \right]^{\sf T}.
\end{equation}
In the simulations, we used
$$P=diag\left(5 \times 10^{-11},0.01,0.005,0.0075, 5 \times 10^{-4}\right),$$
$r_{\tt min}=6628$ km, $U_{\tt max}=10^{-3}$ km/s$^2$,
$e_{\tt min}=10^{-3}$,
 $\epsilon_1=25$, 
 $\epsilon_2=5 \times 10^{-4}$.  
 Aggressive tuning of $P$ was deliberately chosen to induce the activation of constraints and to demonstrate that they are being enforced with the proposed approach based on adding barrier functions to shape the Lyapunov function. Simulation results are shown in Figure~\ref{fig:1} which indicates that constraints are enforced.

\begin{figure}[htpb]
    \centering
\subfloat[]{
    \includegraphics[width=7cm]{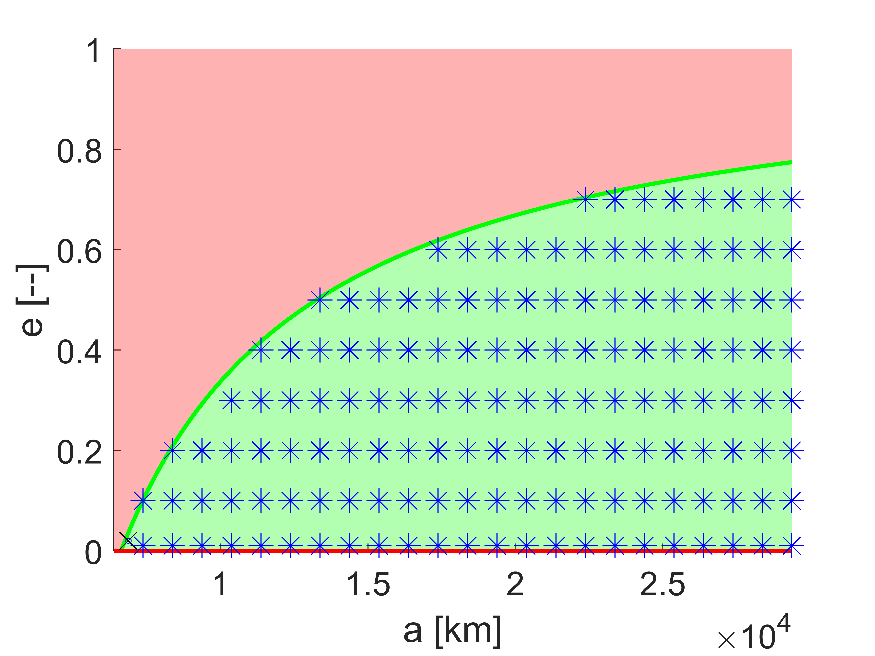}
}\qquad
\subfloat[]{
    \includegraphics[width=7cm]{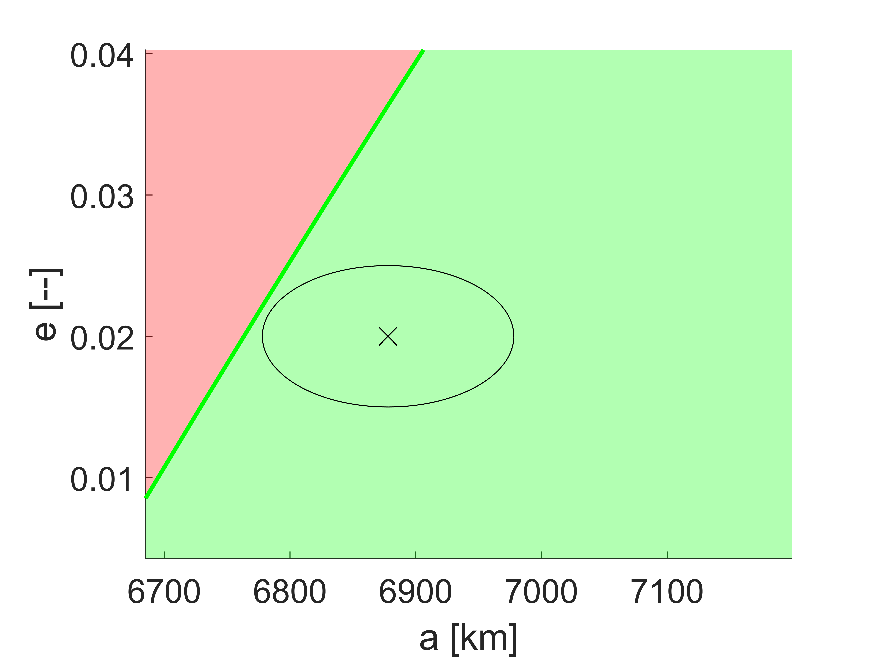}
}
    \caption{(a) Initial conditions (blue ``$\ast$'') which result in closed-loop trajectories convergent to the desired state, $X_{\tt des}$ (black ``X'') within the simulation horizon of 40 hours. (b) Terminal ellipsoid (black) around $X_{\tt des}$ (black ``X''). }
    \label{fig:grid}
\end{figure}

When using barrier functions to avoid constraint violations, theoretical closed-loop stability/ convergence (i.e., liveness)  guarantees
may not be available, however, convergence can be tested with simulations.
For instance, 
we have simulated the closed-loop response to initial conditions corresponding to various values  of the semi-major axis and eccentricity on the chosen mesh.
The convergence was considered successful if the trajectory was able to enter a small ellipsoid around the desired state where constraints are inactive. The mesh prescribed values of $a$ every $1000$ km and values of $e$ every $0.1$, and we have only considered initial conditions instantaneously feasible under constraints.  Figure~\ref{fig:grid} illustrates that given sufficient time, all of the trajectories corresponding to our initial conditions converged into the desired terminal ellipsoid.


\begin{figure}[htpb]
    \centering
    \includegraphics[width = 8 cm]{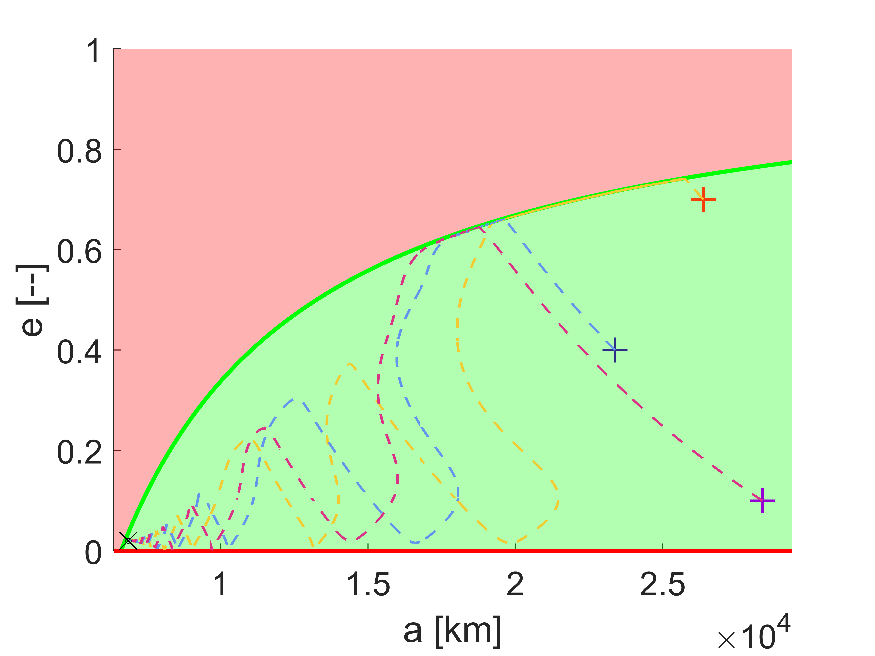}
    \caption{Trajectories in the $a$-$e$ plane for three initial conditions marked by ``+''.}
    \label{fig:examples}
\end{figure}


\section{Integrating Reference Governor into Barrier Function based Lyapunov design to ensure convergence to target orbit (liveness)}


The barrier functions may  create spurious equilibria and deadlocks
away from the target equilibrium,
at which the repulsive action of the barrier function terms to avoid constraint violation is balanced 
by the attracting action of the nominal feedback law.  Closed-loop trajectories could also, in principle,  exhibit limit cycles.  While our simulations did not reveal such issues for the closed-loop trajectories under our specific controller tuning, our use of a finite mesh may have missed problematic initial conditions and target states.  
Hence we propose to use the reference governor \cite{garone2017reference,garone2022command,gilbert2002nonlinear,bemporad1998reference}  to ensure that the closed-loop trajectories 
converge to the target orbit, i.e., $X(t) \to X_{\tt des}$ as $t \to \infty$.

The reference governor operates by temporarily 
replacing the target state, $X_{\tt des}$, by a virtual target state, 
$\tilde{X}_{\tt des}$, that evolves according to
\begin{equation}\label{equ:dec0}
\tilde{X}_{\tt des} (t_n) = \tilde{X}_{\tt des} (t_{n-1})
+ \kappa(t_n) \big(X_{\tt des} -  \tilde{X}_{\tt des} (t_{n-1})\big),~~ 0 \leq \kappa(t_n) \leq 1,~n>1.
\end{equation}
Here $\{t_n\}$ designates a sequence of the time instants at which the target state adjustments are made and $\kappa(t_n)$ is selected by the reference governor algorithm as further delineated below.   

Traditionally, the reference governor is applied to enforce all the constraints, and it can, in principle, be used\cite{semeraro2022reference} in such a capacity for our orbital transfer problem.  
However, in this paper the constraints are enforced through the use of barrier functions
and control saturation; hence, here we only consider the minimal use of the reference governor as a {\em convergence governor}.  

More specifically, our reference governor ensures that the predicted closed-loop trajectory at the end of a chosen prediction horizon, i.e., at $t=t_n+t_{\tt hor}$, $t_{\tt hor}>0$, is guaranteed to enter into a small positively-invariant terminal set around $\tilde{X}_{\tt des}(t_n)$ in which the barrier functions take zero values.


To simplify the developments, we consider the use of the controller
(\ref{equ:controllaw2})-(\ref{equ:controllaw1}) with constant-in-time weights $q_1$ and $q_2$.   
Let 
\begin{equation}\label{equ:dec00}
\hat{X} \left( t_n + t_{\tt hor}; t_n, X(t_n),\theta(t_n),\tilde{X}_{\tt des} \right),
\end{equation}
denote the predicted state of the closed-loop system at the time instant $t_n+t_{\tt hor}>t_n$ where $t_{\tt hor}>0$ is the chosen prediction horizon.  In (\ref{equ:dec00}) we assume that the target state is set to $\tilde{X}_{\tt des}$, the  state
at the time instant $t_n$ is $X(t_n)$,
and the true anomaly at the time instant $t_n$ is  $\theta(t_n)$.

Consider the sublevel sets of the nominal Lyapunov function without the barrier function terms (\ref{equ:V_0}) defined as
$$ Q(\tilde{X}_{\tt des})=\left\{X \in \mathbb{R}^5:~ V_0(P, X,\tilde{X}_{\tt des}) \leq c_0(\tilde{X}_{\tt des}) \right\},$$
where $c_0(\tilde{X}_{\tt des})>0$.  
Suppose that the prediction horizon $t_{\tt hor}>0$
is chosen sufficiently large and the values of  $c_0(\tilde{X}_{\tt des})>0$ are chosen sufficiently small so that the following properties hold:
\begin{description}

\item[(A1)] {\bf Barrier functions are inactive in $Q(\tilde{X}_{\tt des})$ for all $\tilde{X}_{\tt des}$}: 
If $X \in Q(\tilde{X}_{\tt des})$, then
$B_1(X)=B_2(X)=0$;

\item[(A2)] {\bf Small target state adjustments are feasible once in $Q(\tilde{X}_{\tt des})$}: There exists $\epsilon>0$ such that if, for any $t$, $X(t) \in Q(\tilde{X}_{\tt des})$, then \begin{equation}\label{equ:condition3}
\hat{X} \left(t + t_{\tt hor};t,  X(t),\theta,\tilde{X}_{\tt des}+\bar{X} \right) \in Q(\tilde{X}_{\tt des}+\bar{X}) , 
\end{equation}
for all $\theta \in [0,2 \pi)$ and all $\bar{X} \in \mathbb{R}^5$ such that $\|\bar{X}\| \leq \epsilon$;

\item[(A3)] {\bf Control constraints are inactive in $Q({X}_{\tt des})$}: For all $X \in Q(X_{\tt des}) $ and all $\theta \in [0, 2\pi)$, 
$\|U_{\tt nom}(X,X_{\tt des},\theta)\| \leq U_{\tt max}$.

\end{description}

The reference governor algorithm for selecting $\kappa(t_n)$ at the time instant $t_n$ involves maximizing $\kappa(t_n)$, subject to (\ref{equ:dec0})
and subject to
\begin{equation}\label{equ:dec1}
\hat{X} \left( t_n + t_{\tt hor}; t_n, X(t_n),\theta(t_n),\tilde{X}_{\tt des}(t_n) \right)  \in
Q \left(\tilde{X}_{\tt des}(t_n) \right).
\end{equation}
Since this  optimization problem involves maximizing only a scalar parameter $\kappa(t_n)$,  it can be solved using simple bisections.  The state 
(\ref{equ:dec1}) is predicted using forward simulations performed online.  

Assuming that the initial state $X(t_0)$ satisfies the constraints, there is always a feasible choice $\tilde{X}_{\tt des}(t_0)=X(t_0)$ at the time instant $t_0$. 
Note also that  
for $X(t) \in Q(\tilde{X}_{\tt des}(t_n)),$ as a consequence of (A1), $V(x(t))=V_0(P,X(t),\tilde{X}_{\tt des}(t_n))$ and hence
$\frac{dV_0(P,X(t), X_{\tt des}(t_n))}{dt} = \frac{dV(X(t))}{dt} \leq 0$
even if the control input is saturated.  Thus the set $Q(\tilde{X}_{\tt des}(t_n)) $ is positively-invariant.
The condition (\ref{equ:dec1}) implies that if $\tilde{X}_{\tt des}$ is not adjusted and remains equal to $\tilde{X}_{\tt des}(t_n)$ ,  the closed-loop trajectory will enter the terminal set  $Q(\tilde{X}_{\tt des}(t_n))$ and is guaranteed to remain in this set forever after. 
Then at the time instant $t_n$, $n>0$, the previous virtual target state, $\tilde{X}_{\tt des}(t_{n-1})$, is feasible due to the condition (\ref{equ:dec1}) enforced at the time instant $t_{n-1}$.  Thus $\kappa(t_n)=0$ is always feasible.

 Per (A2), once the trajectory enters $Q(\tilde{X}_{\tt des}(t_n))$, sufficiently small in magnitude adjustment of $\tilde{X}_{\tt des}$  towards $X_{\tt des}$ become feasible.  The finite-time convergence of 
$\tilde{X}_{\tt des}(t_n)$ to $X_{\tt des}$
then follows under an additional technical modification\cite{gilbert2002nonlinear} (typically not required in practice) that very small target state adjustments are rejected, i.e.,
if  $\|\tilde{X}_{\tt des}(t_n)-\tilde{X}_{\tt des}(t_{n-1})\|< \delta$, where $\delta<\epsilon$ is small, then $\tilde{X}_{\tt des}(t_n)=\tilde{X}_{\tt des}(t_{n-1})$.  
The finite time convergence of $\tilde{X}_{\tt des}(t_n)$ to $X_{\tt des}$
then implies $X(t) \in Q(X_{\tt des})$ for all $t$ sufficiently large.

By (A3), we now conclude that the asymptotic convergence $X(t) \to {X}_{\tt des}$ as $t \to \infty$  follows under the same conditions as for the closed-loop stability with the nominal unconstrained control law, 
\begin{equation}\label{equ:unc}
U=-G^{\sf T} (X,\theta)
P (X-X_{\tt des})
\end{equation}
without the barrier function terms and control saturation.
A sufficient condition for the latter\cite{semeraro2022reference}  is the persistence of excitation by $G(X(t),\theta(t))$. This persistence of 
excitation condition appears to hold for all of our simulated maneuvers.

\begin{figure}[htpb]
    \centering
    \includegraphics[width = 8 cm]{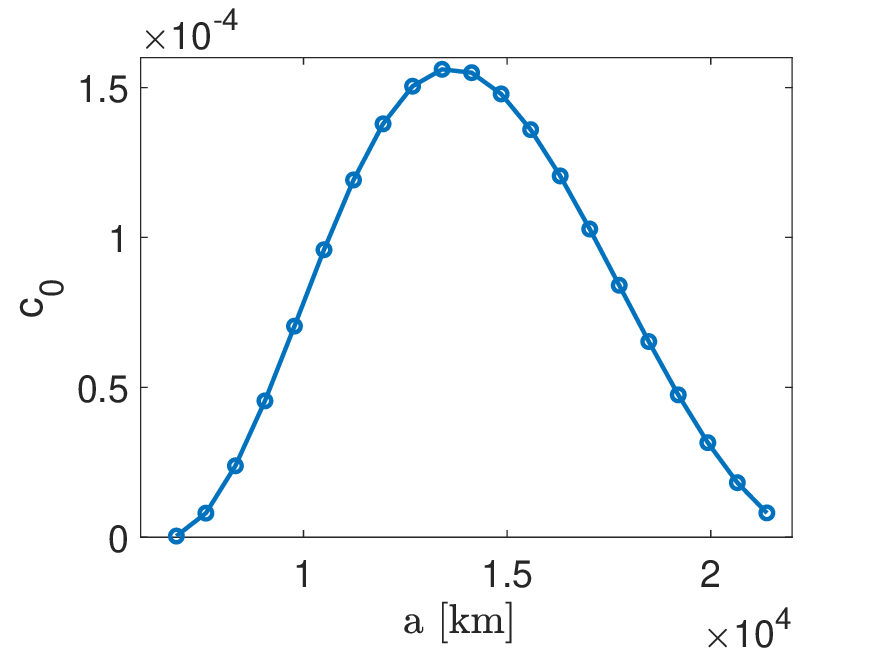}
    \caption{$c_0(\tilde{X}_{\tt des})$ versus $a$. Note that at $a= 6878$ km (leftmost point), $c_0(\tilde{X}_{\tt des})=3.97023 \times 10^{-7} >0$. }
    \label{fig:c0}
\end{figure}

For $X(0)$, $\theta(0)$ and $X_{\tt des}$ defined by (\ref{equ:X0}), (\ref{equ:Xdes}) 
and the same constraints as  for the simulation results in Figure~\ref{fig:1}, the values of $c_0$ versus $a$ are plotted in Figure~\ref{fig:c0}.  These values were computed by selecting virtual target states $\tilde{X}_{\tt des}$ along the line segment connecting $X(0)$  and $X_{\tt des}$ and maximizing numerically
$f_1(X)=r_{\tt min}+\epsilon_1 - a(1 - e)$ and
 $f_2(X)=e_{\tt min}+\epsilon_2 -e $ over a sublevel set of $V_0$ as the size of the sublevel set, determined by $c_0(\tilde{X}_{\tt des})$, was varied. This resulted in the largest $c_0(\tilde{X}_{\tt des})$ for which both the maximum of $f_1(X)$ and the maximum of $f_2(X)$ were non-positive so that (A1) is satisfied. As $Q(X_{\tt des})$ is small-sized, (A3) holds. Condition (A2) has been checked with simulations.

\begin{figure}[htpb]
    \centering
\subfloat[]{
    \includegraphics[width=4cm]{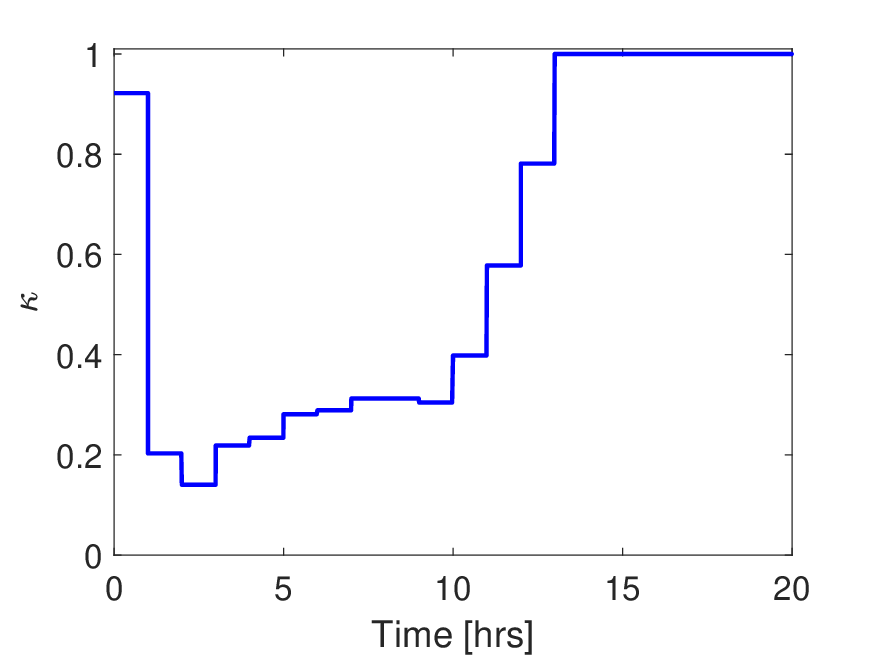}
}\qquad
\subfloat[]{
    \includegraphics[width=4cm]{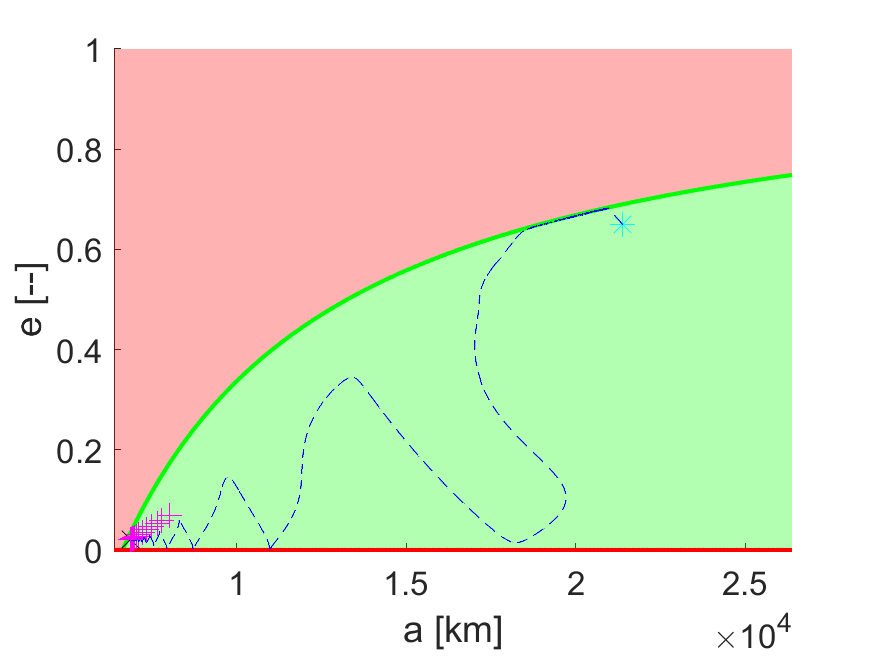}
}
\subfloat[]{
    \includegraphics[width=4cm]{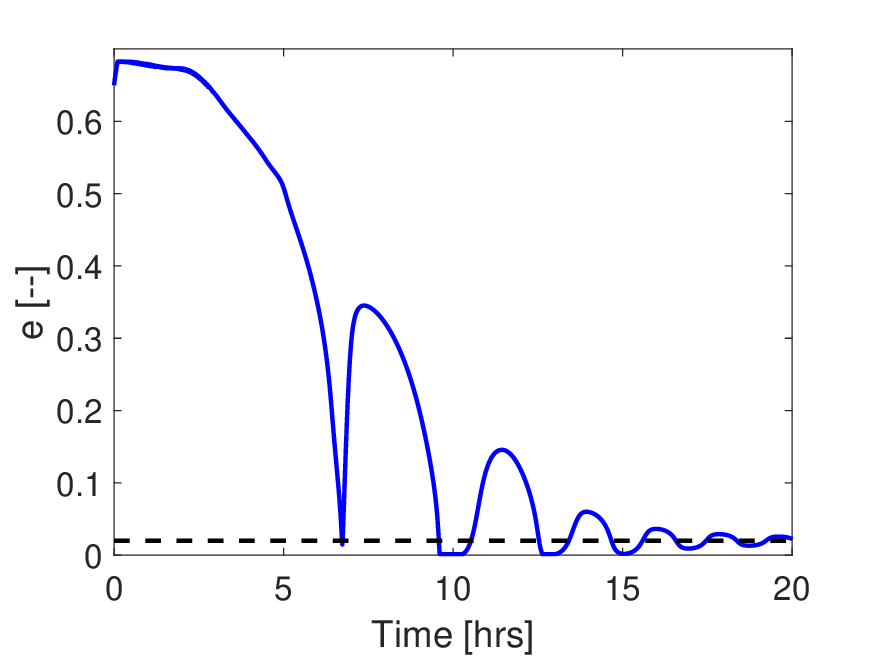}
}
\\
\subfloat[]{
    \includegraphics[width=4cm]{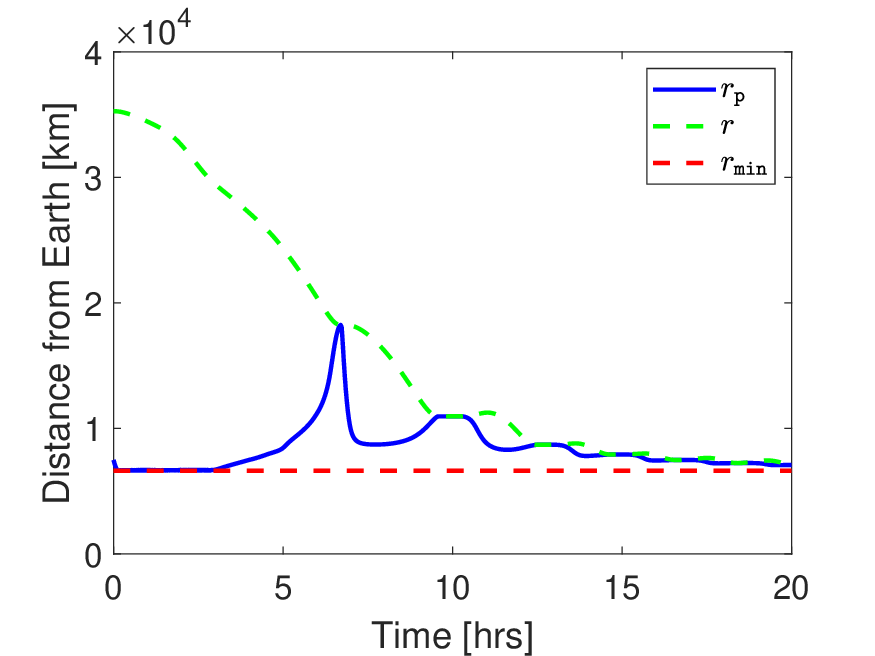}
}
\subfloat[]{
    \includegraphics[width=4cm]{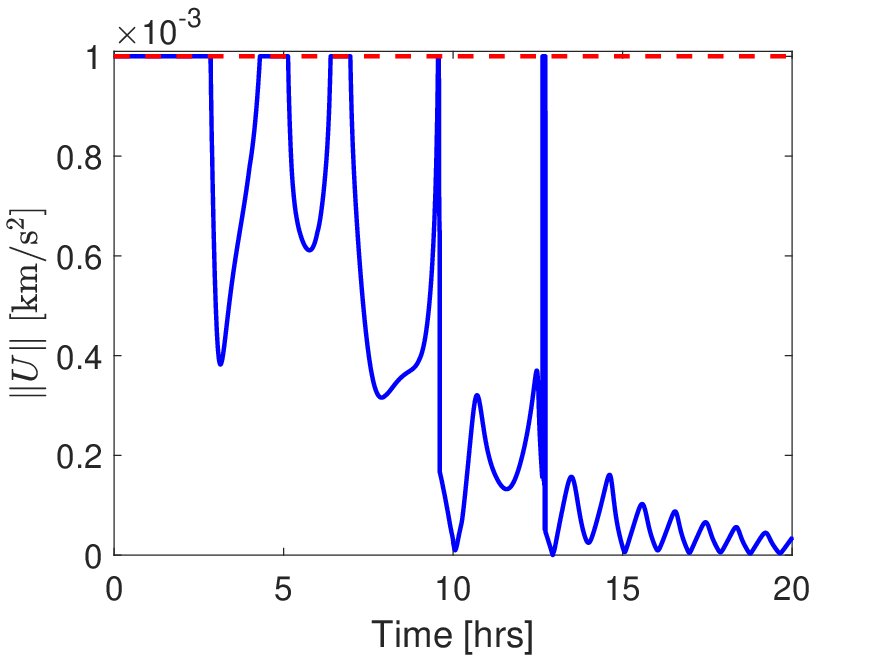}
}
    \caption{Response with reference governor for $t_{\tt hor}=15$ hours: (a) Time history of $\kappa(t)$; (b) 
     Trajectory in the $a$-$e$ plane. The magenta $+$ indicates the virtual target; (c) Time history of eccentricity; (d) The time histories of $r_{\tt p}$, $r$ and $r_{\tt min}$; (e)  The time histories of $\|U \|$ and  $U_{\tt max}$. }
    \label{fig:RG}
\end{figure}

Figure~\ref{fig:RG} shows the responses with the reference governor where the prediction horizon is $t_{\tt hor}=15$
hours.  Note that the response is essentially unchanged by the addition of the reference governor; however, the reference governor enforces convergence.  

\begin{figure}[htpb]
    \centering
\subfloat[]{
    \includegraphics[width=4cm]{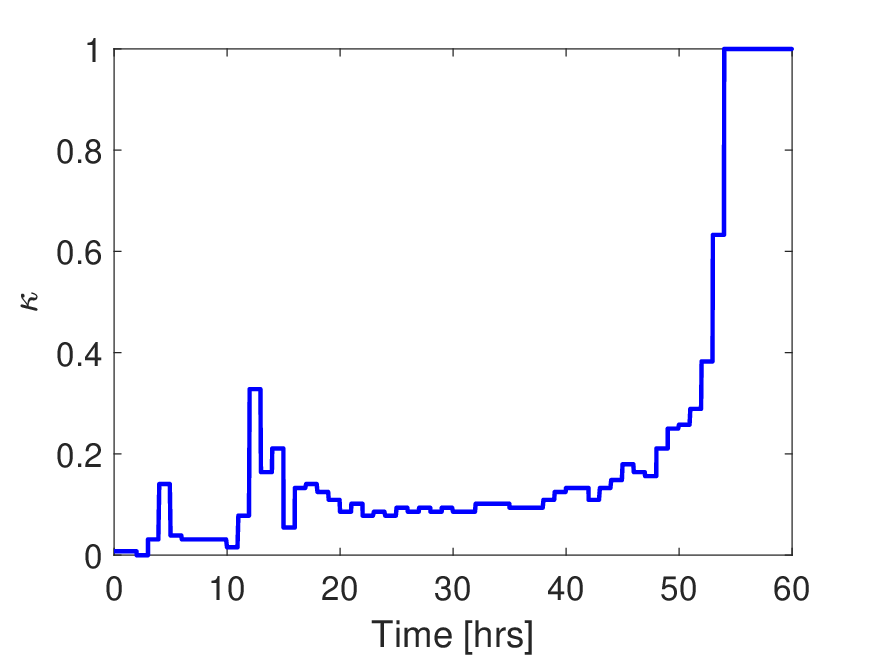}
}\qquad
\subfloat[]{
    \includegraphics[width=4cm]{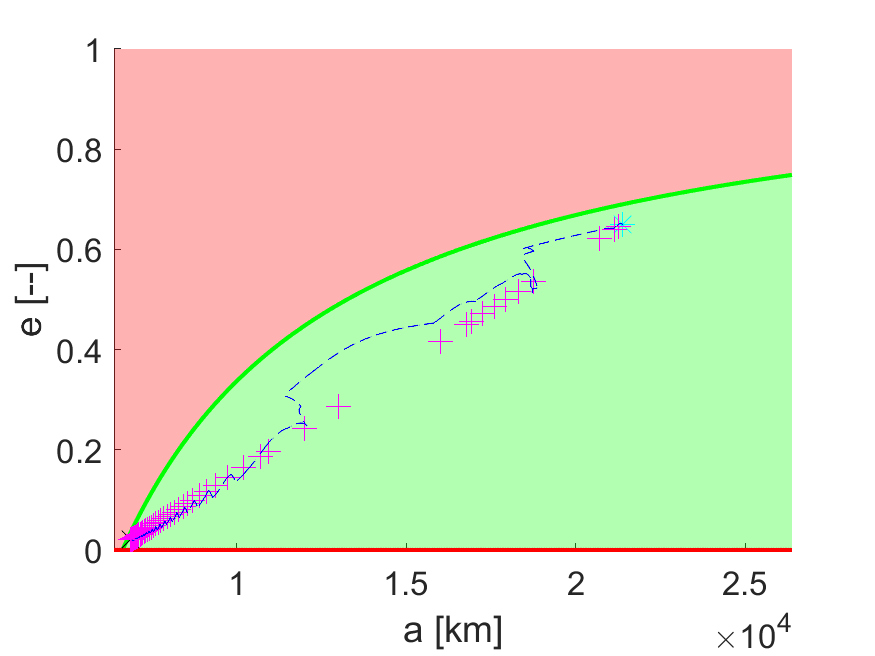}
}
\subfloat[]{
    \includegraphics[width=4cm]{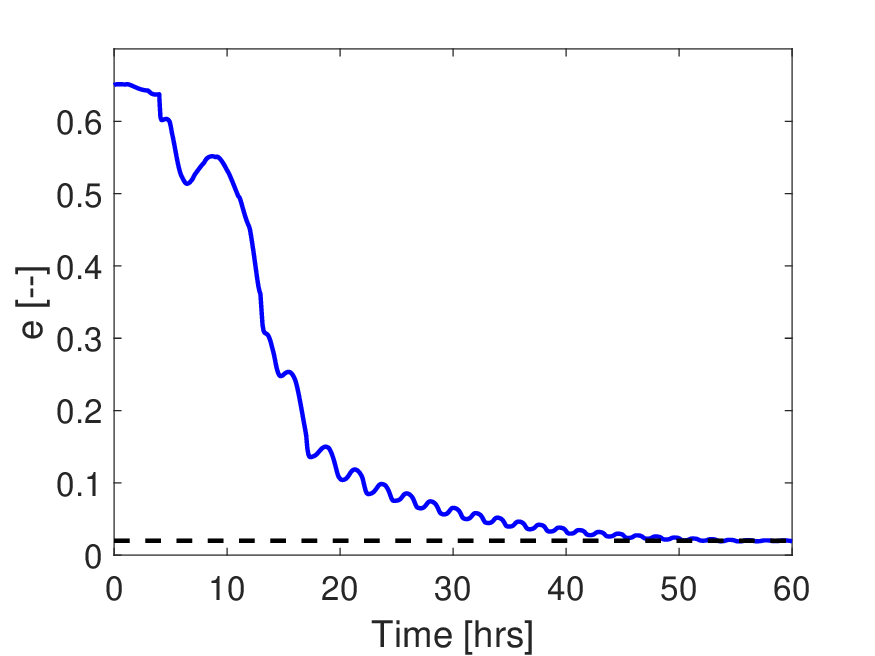}
}
\\
\subfloat[]{
    \includegraphics[width=4cm]{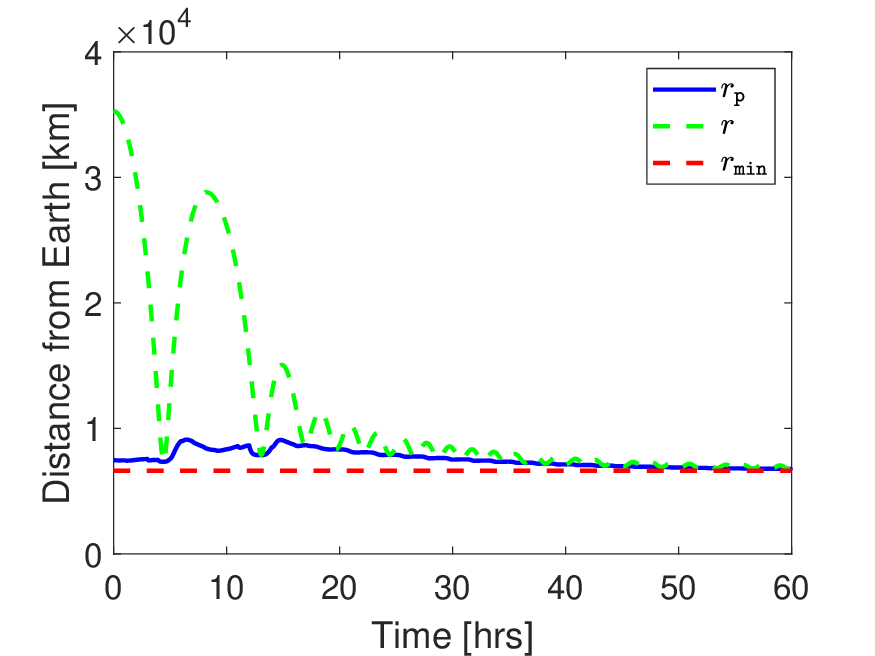}
}
\subfloat[]{
    \includegraphics[width=4cm]{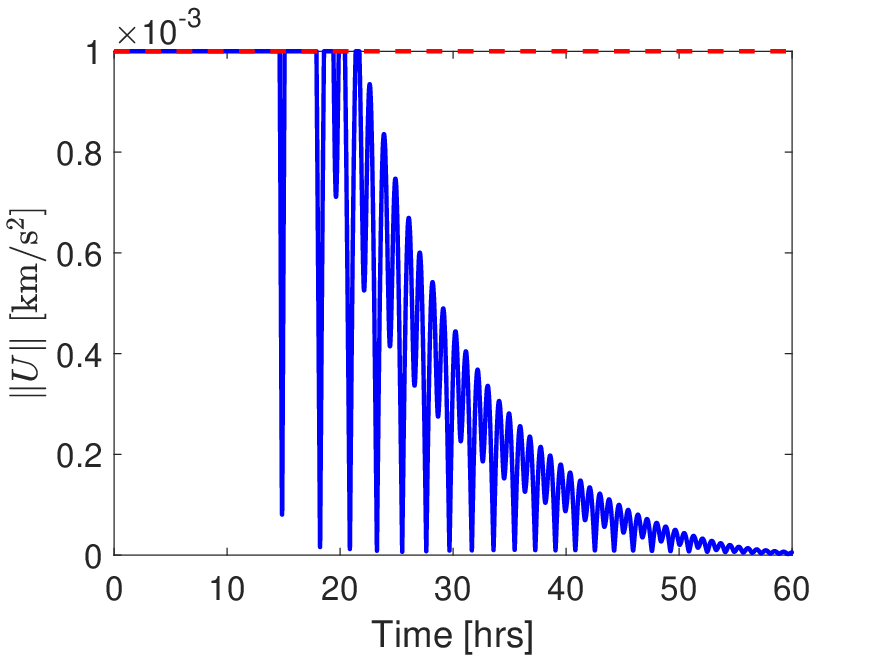}
}
    \caption{Response with reference governor and $t_{\tt hor}=4$ hours: (a) Time history of $\kappa(t)$; (b) 
     Trajectory in the $a$-$e$ plane. The magenta $+$ indicates the virtual target; (c) Time history of eccentricity; (d) The time histories of $r_{\tt p}$, $r$ and $r_{\tt min}$; (e)  The time histories of $\|U \|$ and  $U_{\tt max}$. }
    \label{fig:RG_shorter}
\end{figure}

The use of shorter prediction horizon, while beneficial in terms of reducing the computational load since the forward prediction is performed over a shorter horizon, at the same time reduces the response speed. See Figure~\ref{fig:RG_shorter} for the response with the prediction horizon, $t_{\tt hor}=4$ hours.  This reduction is expected as the terminal constraint must be satisfied over a much shorter time duration reducing the changes in $\tilde{X}_{\tt des}$ that the reference governor can take at each step.


\section{Conclusions}
Lyapunov-based feedback controllers can be developed for spacecraft orbital maneuvering which satisfy state and control constraints through the application of the barrier functions and control input saturation. 
Reference governor can be used as a convergence governor to avoid spurious equilibria and limit cycles that can occur with barrier function based methods and to ensure convergence to the target orbit.
The results are dependent on and exploit the drift-free form of the spacecraft equations of motion when expressed in terms of Gauss Variational Equations.

\bibliographystyle{AAS_publication}   
\bibliography{references}   

\section*{Appendix: Handling more general constraints  on the control input}

We first illustrate the derivation of a Lyapunov controller for the $\infty$-norm constraints on the control input of the form,
$$|U_i(t)| \leq U_{{\tt max},i},\quad i=1,2,3,$$
that can be written concisely as $U\in \mathcal{U}=\mathcal{U}_1 \times \mathcal{U}_2 \times \mathcal{U}_3$, where $ \mathcal{U}_i=[-U_{max,i},U_{max,i}], i=1,2,3. $     

For the moment, we do not consider other constraints or the use of barrier functions to enforce them. We will come back to  the more general case at the end of the Appendix. 

Given the equations (\ref{gve}) written in the form,
$$    \dot X = G(X(t),\theta(t))U(t), $$
we consider  the nominal control law,
$$ U_{\tt nom}(X,X_{\tt des})=-G(X,\theta)^{\sf T} P (X-X_{\tt des})=-G^{\sf T}P E =-G^{\sf T} \tilde{E},
$$
where we use $G$ as a shorthand for $G(X(t),\theta(t))$, 
while $E=(X-X_{\tt des})$
and $\Tilde{E}=PE$.

The saturated control input is given by
$$    U_i = \Pi_{\mathcal{U}_i}[-G_i^{\sf T}P(X-X_{\tt des})] = \left\{ {\begin{array}{*{20}l} {-G_i^{\sf T}\Tilde{E},} \hfill & {{\text{if}} \mid G_i^{\sf T}\Tilde{E}\mid \leq  U_{{\tt max},i} } \hfill \\ {-U_{{\tt max},i}sign\left(G_i^{\sf T}\Tilde{E}\right),} \hfill & {{\text{otherwise,}}} \hfill \\ \end{array} } \right. \\
$$
where $$G^{\sf T} =\left[\begin{array}{c} G_1^{\sf T} \\ G_2^{\sf T} \\ G_3^{\sf T} \end{array} \right].$$
Then defining the Lyapunov function candidate by
$$
    V = \frac{1}{2}\Tilde{E}^{\sf T}P^{-1}\Tilde{E}^{\sf T}. 
$$    
and computing its time derivative along the closed-loop system trajectories we obtain,
\begin{gather*}
    \dot V=\Tilde{E}^{\sf T}P^{-1}PGU=\Tilde{E}^{\sf T}G_1 U_1+\Tilde{E}^{\sf T}G_2 U_2+\Tilde{E}^{\sf T}G_3 U_3 \\ 
    = \sum_{i=1}^3 \Tilde{E}^{\sf T}G_i\left\{ {\begin{array}{*{20}l} {-G_i^{\sf T}\Tilde{E},} \hfill & {{\text{if}} \mid G_i^{\sf T}\Tilde{E}\mid \leq  U_{{\tt max},i} } \hfill \\ {-U_{{\tt max},i}sign\left(G_i^{\sf T}\Tilde{E}\right),} \hfill & {{\text{otherwise}}} \hfill \\ \end{array} } \right.  \\
    = \sum_{i=1}^3 \left\{ {\begin{array}{*{20}l} {-\Tilde{E}^{\sf T}G_iG_i^{\sf T}\Tilde{E},} \hfill & {{\text{if}} \mid G_i^{\sf T}\Tilde{E}\mid \leq  U_{{\tt max},i} } \hfill \\ {-U_{{\tt max},i}\mid G_i^{\sf T}\Tilde{E}\mid,} \hfill & {{\text{otherwise}}} \hfill \\ \end{array} } \right. \leq 0.
\end{gather*}

\begin{figure}[htpb]
\centering
\subfloat[]{
    \includegraphics[width=4cm]{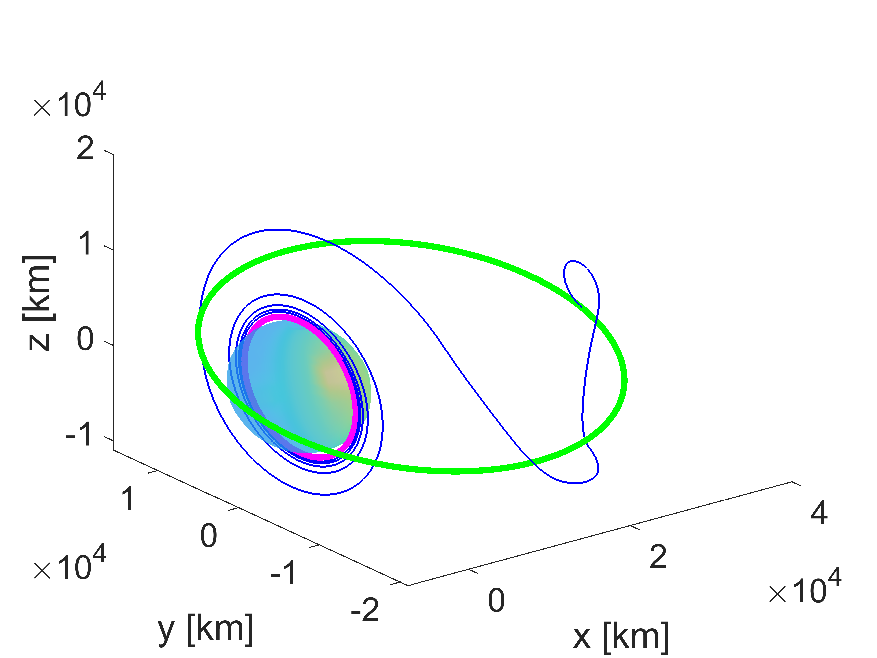}
}\qquad
\subfloat[]{
    \includegraphics[width=4cm]{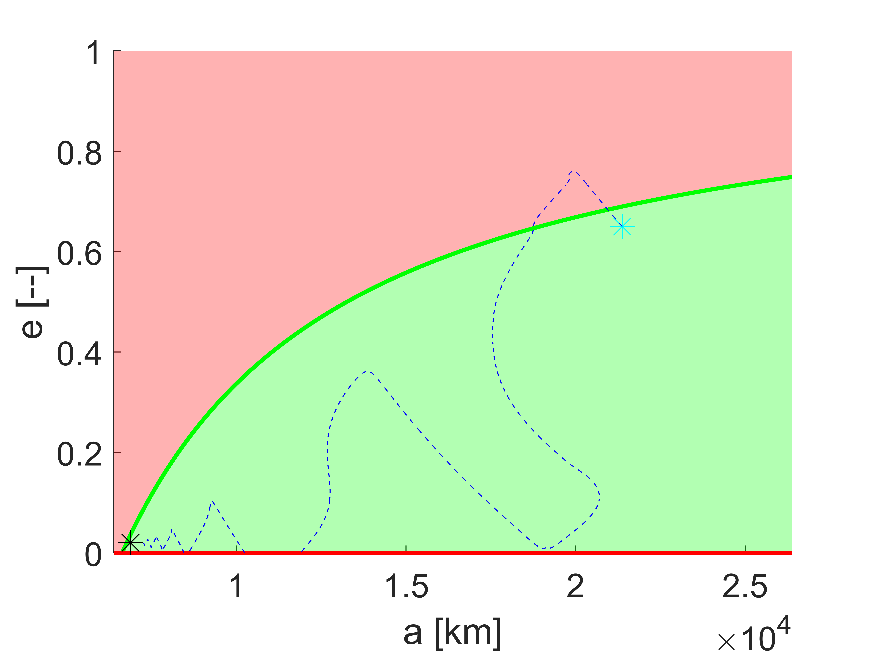}
}\qquad
\subfloat[]{
    \includegraphics[width=4cm]{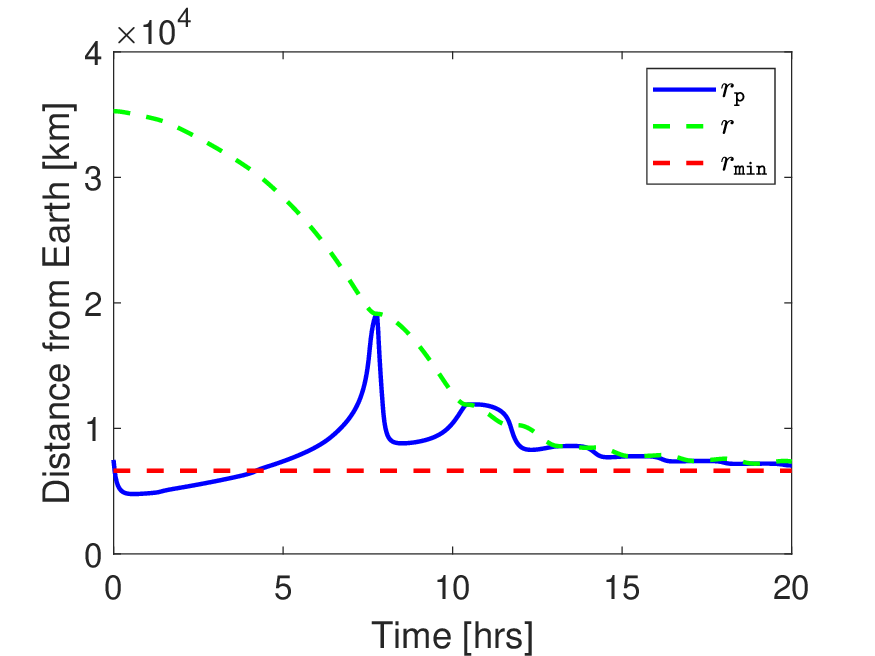}
}\vskip\baselineskip
\subfloat[]{
    \includegraphics[width=4cm]{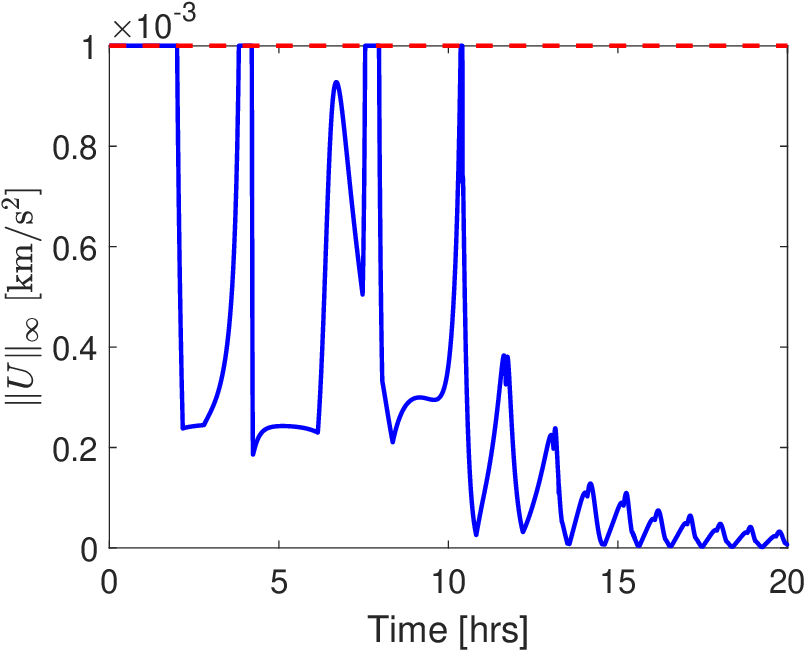}
}\qquad
\subfloat[]{
    \includegraphics[width=4cm]{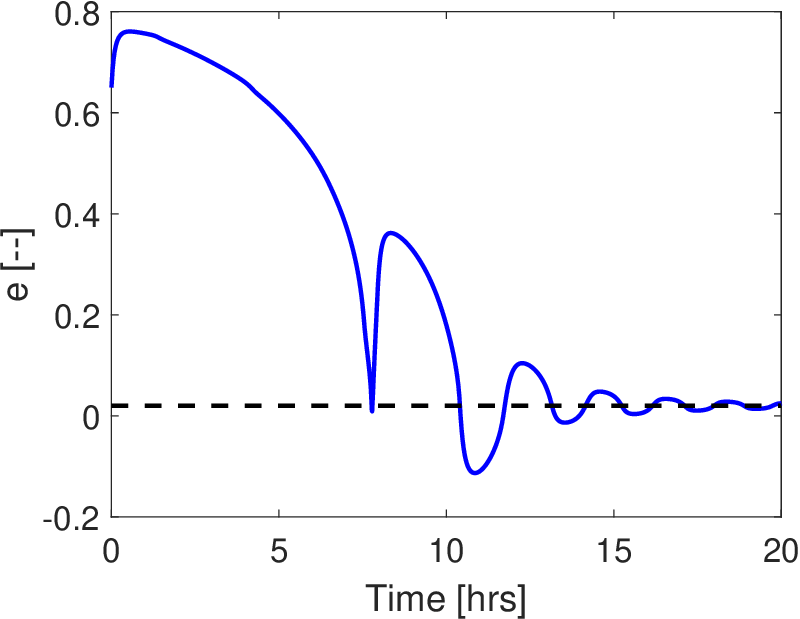}
}\qquad
\subfloat[]{
    \includegraphics[width=4cm]{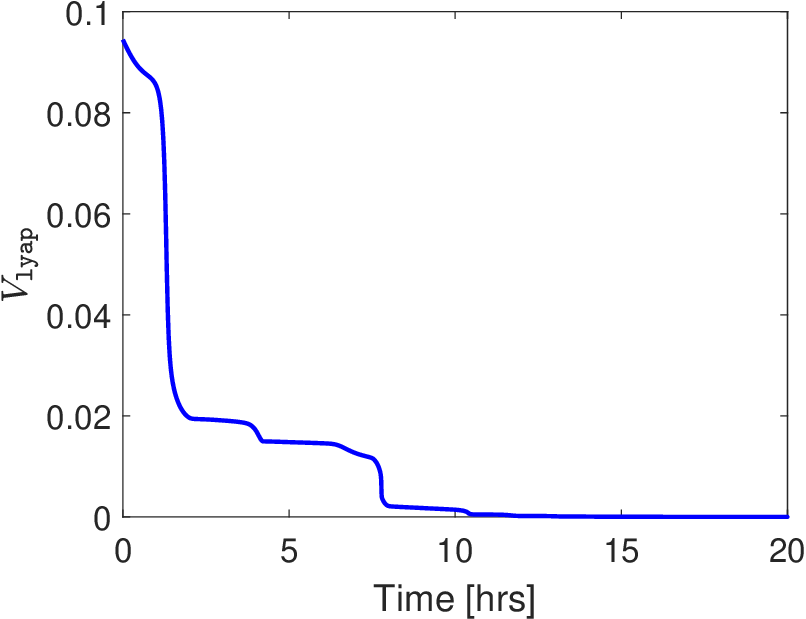}
}\qquad
\subfloat[]{
    \includegraphics[width=4cm]{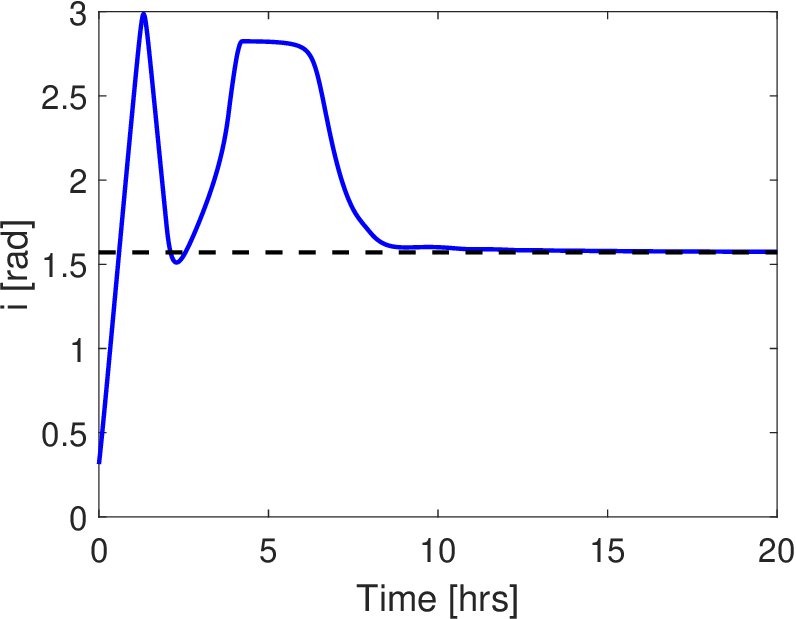}
}\qquad
\subfloat[]{
    \includegraphics[width=4cm]{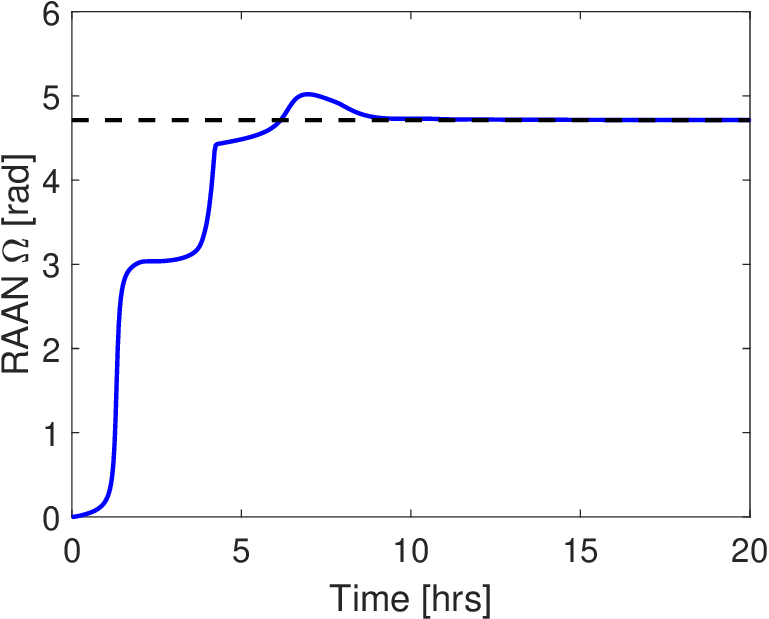}
}\qquad    
\subfloat[]{
    \includegraphics[width=4cm]{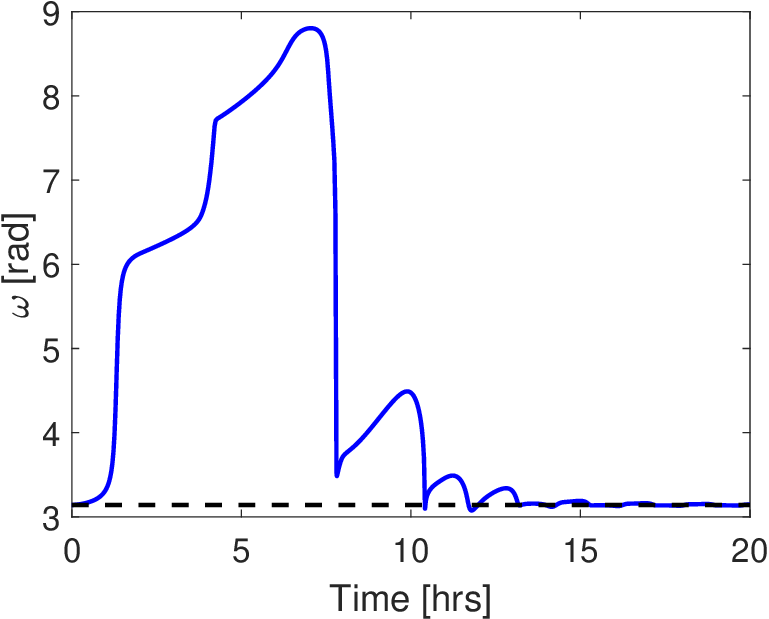}
}
\caption{Orbital transfer from a higher orbit to a lower orbit with only $\infty$-norm control constraint: (a) Three dimensional trajectory; (b) Trajectory on $a$-$e$ plane (dashed) with the region allowed by constraints (\ref{equ:maincon1}) and (\ref{equ:maincon3}) shown in green; (c) The time histories of $r_{\tt p}$, $r$ and $r_{\tt min}$; (d) The time histories of $\|U\|_\infty$
and $U_{\tt max}=U_{{\tt max},i}$, $i=1,2,3$; (e) The time histories of eccentricity $e$ (solid) and $e_{\tt des}$ (dashed); (f) The time history of the Lyapunov function, $V(X(t)$;
(g) The time histories of inclination $i$ (solid) and $i_{\tt des}$ (dashed);
(h) The time histories of RAAN $\Omega$ (solid) and $\Omega_{\tt des}$ (dashed);
(i) The time histories of argument of periapsis $\omega$ (solid) and $\omega_{\tt des}$ (dashed).}\label{fig:3}
\end{figure}

We illustrate the closed-loop response in Figure~\ref{fig:3}.  Note that $X(t) \to X_{\tt des}
$ as $t \to \infty$ and 
control constraints are enforced.  State constraints (\ref{equ:maincon1}) and (\ref{equ:maincon3})  are violated as only control constraints are imposed in this case.

Consider now a more general case when the barrier functions are used
to enforce state constraints, and control constraints have the form, $$U(t) \in \mathcal{U},$$
where $\mathcal{U}$ is a convex, compact set with $0 \in \mathcal{U}$.  The control law in this case can be defined using the minimum $2$-norm projection onto the set $\mathcal{U}$ as
$$    U = \Pi_{\mathcal{U}}[U_{\tt nom}],\quad U_{\tt nom} =  
-\left[(X-X_{\tt des})^{\sf T} PG(X,\theta) + C^{\sf T}(X)G(X,\theta)\right]^{\sf T}.
$$ 
Since $U$ is the minimum norm projection of $U_{\tt nom}$ onto $\mathcal{U}$ and $0 \in \mathcal{U}$, which is convex, the necessary conditions 
for optimality of $U$ in the problem of minimizing $f(u)=\|u-U_{\tt nom}\|^2$ subject to $u \in \mathcal{U}$
imply
$$  (\nabla_u f(U))^{\sf T}(v-U)  \leq 0, \quad {\mbox{for any $v \in \mathcal{U}$}} $$
and with $v=0$ that
$$ (U_{\tt nom} - U)^{\sf T} (0-U) \leq 0 \Rightarrow
-U_{\tt nom}^{\sf T} U \leq - U^{\sf T} U.
$$
Thus, with $V$ given by (\ref{equ:main2}), it follows that
its time rate of change of $V(X(t))$ along the trajectories of the closed-loop system satisfies
\begin{gather*}
    \frac{d}{dt} V(X(t)) = \left((X(t)-X_{des})^{\sf T} PG(X(t),\theta(t)) + C^{\sf T} (X(t))G(X(t),\theta(t)) \right) U(t) \\
    =-U_{\tt nom}^{\sf T}(t) U(t)  \leq -U^{\sf T}(t) U(t) \leq 0.
    \end{gather*}
Thus the sublevel sets of $V$ remain invariant, and hence the constraints can be enforced using exactly the same procedure as in the case of the $2$-norm constraints (\ref{equ:maincon2}) on the control input.

\vskip 2em

\end{document}